\documentclass[12pt,draft]{amsart}
\usepackage{amsmath}
\usepackage{amssymb}
\usepackage{amsthm}
\usepackage{bm}
\numberwithin{equation}{section}
\newtheorem{theorem}{Theorem}[section]

\newtheorem{proposition}[theorem]{Proposition}

\allowdisplaybreaks
\begin{document}
\title[Simple proof of the functional relation for the Tornheim double zeta]{Simple proof of the functional relation for the Lerch type Tornheim double zeta function}
\author[T.~Nakamura]{Takashi Nakamura}
\address{Department of Mathematics Faculty of Science and Technology Tokyo University of Science Noda, CHIBA 278-8510 JAPAN}
\email{nakamura\_takashi@ma.noda.tus.ac.jp}
\begin{abstract}
In this paper, we give a simple proof of the functional relation for the Lerch type Tornheim double zeta function. By using it, we obtain simple proofs of some explicit evaluation formulas for double $L$-values. 
\end{abstract}
\maketitle

\section{Introduction and main results}
We define the Lerch type Tornheim double zeta function by
\begin{equation}
T (s,t,u \,; x,y) := \lim_{R \to \infty} \sum_{m,n=1}^{m+n=R} \frac{e^{2 \pi i mx} e^{ 2 \pi i ny}}{m^s n^t (m+n)^u},
\label{eq:defT}
\end{equation}
where $0 \le x,y \le 1$, $\Re (s+t) >1$, $\Re (t+u) > 1$ and $\Re (s+t+u) > 2$. This function is continued meromorphically by \cite[Theorem 2.1]{NakamuraA}. Let $k \in {\mathbb{N}} \cup \{ 0 \}$. The function $T (s, t, u \,; x,y)$ can be continued meromorphically to ${\mathbb{C}}^3$, and all of its singularities are located on the subsets of ${\mathbb{C}}^3$ defined by the following equations;
\begin{equation*}
\begin{split}
t = 1 - k \qquad &{\rm{if}} \quad x \not \equiv 1, \,\, y  \equiv 1 \mod 1, \\
s = 1 - k \qquad &{\rm{if}} \quad x \equiv 1, \,\, y \not \equiv 1 \mod 1, \\
{\mbox{no singularity}} \qquad &{\rm{if}} \quad x \not \equiv 1, \,\, y \not \equiv 1 \mod 1.
\end{split}
\end{equation*}

We write $T(s, t, u) := T (s, t, u \,; 1,1)$ and call this function the Tornheim double zeta function. The values $T(a,b,c)$ for $a,b,c \in{\mathbb{N}}$ were investigated by Tornheim in 1950 and later by Mordell in 1958, and some explicit formulas for them were obtained. Subbarao and Sitaramachandrarao, Huard, Williams and Zhang, and Tsumura researched the explicit formulas for $T(a,b,c)$ for $a,b,c \in{\mathbb{N}}$. The value $T (0, a, b \,; x,y)$ and their multiple sum versions have been already defined in Arakawa and Kaneko \cite{AraKa} for the case $x,y \in {\mathbb{Q}}$ as special cases of their multiple $L$-values. 

As a three-variable function, Matsumoto continued $T(s, t, u)$ meromorphically to the whole ${\mathbb{C}}^3$ plane in \cite[Theorem 1]{Ma}. Tsumura \cite[Theorem 4.5]{Tsumurap2}, afterwards Nakamura \cite[Theorem 1.2]{Na} found functional relations for the Tornheim double zeta function. Moreover, generalizations of the functional relations are proved by Matsumoto and Tsumura \cite{MaTsu}, and Nakamura \cite{NakamuraA}. 

In this paper, we show the following functional relation. This functional relation is essentially the same as \cite[Theorem 3.1]{NakamuraT}. Therefore we can obtain the all results in \cite{NakamuraT} by this formula. Zhou gave a simple proof of \cite[Theorem 1.2]{Na} in \cite{Zhou2}. Recently, Li gave the proof similar to Zhou's one in \cite{Li}, independently. By modifying their methods, we can prove the following theorem. 
\begin{theorem}
For $0 < x \ne y < 1$, $a,b \in {\mathbb{N}}$ and $s \in {\mathbb{C}}$, we have
\begin{equation}
\begin{split}
&T (a,b,s \,; x,y) + (-1)^b T(b,s,a \,; x-y,x) + (-1)^a T(s,a,b \,; y, y-x) \\
&=\sum_{j=1}^{a} \binom{a+b-j-1}{a-j} \zeta (a+b+s-j \, ; y ) \bigl( \zeta (j \, ; x-y) + (-1)^j \zeta (j \, ; y-x) \bigr) \\ 
&\quad + \sum_{j=1}^{b} \binom{a+b-j-1}{b-j} \zeta (a+b+s-j \, ; x ) \bigl( \zeta (j \, ; y-x) + (-1)^j  \zeta (j \, ; x-y )\bigr) \\
&\quad - \binom{a+b-1}{a}  \zeta (a+b+s \, ; y) - \binom{a+b-1}{b}  \zeta (a+b+s \, ; x) .
\label{eq:th1}
\end{split}
\end{equation}
\label{th:1}
\end{theorem}

Taking $x \to y$ in the above formula, we have \cite[(3.1)]{NakamuraT} since
$$
\lim_{x \to y} \bigl( \zeta (a+b+s-1 \, ; y ) - \zeta (a+b+s-1 \, ; x ) \bigr) 
\bigl( \zeta (1 \, ; x-y) - \zeta (1 \, ; y-x) \bigr) = 0 .
$$
Define $K (a,b \,; x,y)$ by the right-hand side of (\ref{eq:th1}) with $s=0$. By using Theorem \ref{th:1}, we obtain the following propositions. It should be noted that $T(0,a,b \,; -y,x-y) = T(0,a,b \,; y, y-x)$ when $(x,y) = (1,1)$, $(1,1/2)$, $(1/2,1)$ or $(1/2, 1/2)$. In this case, the next proposition coincides with \cite[Proposition 2.8]{NakamuraT} or \cite[Proposition 1]{Zhou}.
\begin{proposition}
For any admissible index, we have
\begin{equation}
\begin{split}
&T(0,a,b \,; -y,x-y) -(-1)^{a+b} T(0,a,b \,; y, y-x)
= \\ &(-1)^b \zeta (a \,; x) \zeta (b \,; y) - (-1)^b K (a,b \,; x,y) 
+ \zeta (a \,; -y)  \zeta (b \,; x) - \zeta (a+b \,; x-y)
\label{eq:tasite}
\end{split}
\end{equation}
\label{pro:1}
\end{proposition}

Let $\phi$, $\chi$ and $\psi$ are Dirichlet characters of conductor $h$, $k$, and $q$, respectively.  For any admissible index, we define $L (0,a,b \,; \phi, \chi ,\psi)$ by
\begin{equation}
L (0,a,b \,; \phi, \chi ,\psi) := \lim_{R \to \infty} \sum_{m,n=1}^{m+n=R}
\frac{\phi (m) \chi (n) \psi (m+n)}{n^a (m+n)^b}.
\label{eq:defL}
\end{equation}

Terhune \cite{Terhune} showed that if and $\chi \psi (-1) = (-1)^{a+b+1}$ then $L (0,a,b \,; 1, \chi ,\psi)$ can be expressed as a polynomial in the values of polylogarithms at certain roots of unity, with coefficients in a cyclotomic field. In \cite[Proposition 4.5]{NakamuraT}, the author obtained explicit evaluation formulas for $L (0,a,b \,; 1, \chi ,\psi)$ when $\chi \psi (-1) = (-1)^{a+b+1}$. Proposition \ref{pro:1} and the following proposition give simpler ones. Denote the Gauss sum by $\tau (\chi) := \sum_{l=1}^{k-1} \chi (l) e^{2 \pi i l/k}$.
\begin{proposition}
Define $2U(a,b \,; x,y):= T(0,a,b \,; x,y) -(-1)^{a+b} T(0,a,b \,; -x, -y)$. Let $\phi \chi \psi(-1) = (-1)^{a+b+1}$. For any admissible index, we have
\begin{equation}
\begin{split}
&\tau (\overline{\phi}) \tau (\overline{\chi}) \tau (\overline{\psi}) L(0,a,b \,; \phi, \chi ,\psi) = \\
&\sum_{j=1}^{h-1} \sum_{l=1}^{k-1} \sum_{r=1}^{q-1} \overline{\phi}(j) \overline{\chi}(l) \overline{\psi}(r) U (a,b \,; j/h+r/q, l/k+r/q).
\label{eq:dl1}
\end{split}
\end{equation}
\label{pro:2}
\end{proposition}

\section{Proof of results}
\begin{proof}[Proof of Theorem \ref{th:1}]
First, suppose $\Re (s) >1$. We define $S (a,b,s \,; x,y)$ by 
$$
S (a,b,s \,; x,y) := T (a,b,s \,; x,y) + (-1)^b T(b,s,a \,; x-y,x) + (-1)^a T(s,a,b \,; y, y-x).
$$
It is easy to see that $S (a,b,s \,; x,y) = S (a-1,b,s+1 \,; x,y) + S (a,b-1,s+1 \,; x,y)$ because of $T (a,b,s \,; x,y) = T (a-1,b,s+1 \,; x,y) + T (a,b-1,s+1 \,; x,y)$. Hence we have
\begin{equation*}
\begin{split}
S (a,b,s \,; x,y) = &\sum_{j=1}^{a} \binom{a+b-j-1}{a-j} S (j,0,a+b+s-j  \,; x,y) \\ &+ \sum_{j=1}^{b} \binom{a+b-j-1}{b-j} S (0,j,a+b+s-j  \,; x,y).
\end{split}
\end{equation*}
Now we consider the function  $S (j,0,a+b+s-j  \,; x,y)$ in the above formula. By the definition of $S (a,b,s \,; x,y)$ and the harmonic product formula, we have 
\begin{equation*}
\begin{split}
&S (j,0,a+b+s-j \,; x,y) = 
T (j,0,a+b+s-j \,; x,y) \\ &+ T(0,a+b+s-j,j \,; x-y,x) + (-1)^j T(a+b+s-j,j,0 \,; y, y-x) \\= \, 
&\zeta (a+b+s-j \, ; y ) \bigl( \zeta (j \, ; x-y) + (-1)^j \zeta (j \, ; y-x) \bigr) - \zeta (a+b+s \, ; x ) .
\end{split}
\end{equation*}
By exchanging the parameters $x$ and $y$, we can repeat the same manner for $S (0,j,a+b+s-j  \,; x,y)$.  Therefore we obtain this theorem when $\Re (s) >1$ by the well-known formula $\sum_{j=1}^{a} \binom{a+b-j-1}{a-j} = \binom{a+b-1}{b}$. By the analytic continuation, we obtain Theorem \ref{th:1}.
\end{proof}

\begin{proof}[Proof of Proposition \ref{pro:1}]
By putting $s=0$ in (\ref{eq:th1}), we have
$$
\zeta (a \,; x) \zeta (b \,; y) + (-1)^b T(b,0,a \,; x-y,x) + (-1)^a T(0,a,b \,; y, y-x) = K (a,b \,; x,y)
$$
On the other hand, one has
$$
T(0,a,b \,; -y,x-y) + T(b,0,a \,; x-y,x) + \zeta (a+b \,; x-y) = \zeta (a \,; -y)  \zeta (b \,; x)
$$
by the harmonic product formula. Therefore we have Proposition \ref{pro:1} by removing the term $T(b,0,a \,; x-y,x)$. 
\end{proof}

\begin{proof}[Proof of Proposition \ref{pro:2}]
Recall the well-known formula
$$
\chi (n) = \frac{1}{\tau (\overline{\chi})} \sum_{l=1}^{k-1} \overline{\chi} (l) e^{2 \pi i ln/k} =
\frac{\chi (-1)}{\tau (\overline{\chi})} \sum_{l=1}^{k-1} \overline{\chi} (l) e^{-2 \pi i ln/k}.
$$
By using above formula and $\phi \chi \psi(-1) = (-1)^{a+b+1}$, we have
\begin{equation*}
\begin{split}
&\tau (\overline{\phi}) \tau (\overline{\chi}) \tau (\overline{\psi}) L(0,a,b \,; \phi, \chi ,\psi)  \\
&=\sum_{j=1}^{h-1} \sum_{l=1}^{k-1} \sum_{r=1}^{q-1} \overline{\phi}(j) \overline{\chi}(l) \overline{\psi}(r) T (0,a,b \,; j/h+r/q, l/k+r/q)\\
&= (-1)^{a+b+1} \sum_{j=1}^{h-1} \sum_{l=1}^{k-1} \sum_{r=1}^{q-1} \overline{\phi}(j) \overline{\chi}(l) \overline{\psi}(r) T (0,a,b \,; -j/h-r/q, -l/k-r/q) .
\end{split}
\end{equation*}
Hence we obtain Proposition \ref{pro:2}. 
\end{proof}

 
\end{document}